\documentclass[12pt,a4paper,oneside,onecolumn]{article}
\title{
On a Poisson-algebraic characterization of vector bundles}
\author{Zihindula Mushengezi E.}
\usepackage{etex}
\usepackage[]{latexsym}
\usepackage[T1]{fontenc}
\usepackage[latin1]{inputenc}
\usepackage[english]{babel}
\usepackage{graphics}
\usepackage[]{amssymb}
\usepackage{multicol}
\usepackage{amsmath}
\usepackage{graphicx}
\usepackage{amsfonts}
\usepackage[all]{xy}
\usepackage{pifont}
\usepackage{microtype}
\usepackage{answers}
\usepackage{tabularx}
\usepackage{float}
\usepackage{color}

\usepackage{listings}
\definecolor{darkWhite}{rgb}{0.94,0.94,0.94}
\usepackage{palatino}

\newcount\exenum     \newdimen\numindent
\newcount\qqnum     \newdimen\qqmarge
\newcount\sqnum     \newdimen\sqmarge

\outer\def\bye{%
\vskip 0pt\@endmulticol\@endgroup              
\ifanswer \let\next=\exobye@                   
\else     \let\next=\@exobye                   
\fi\next}

\usepackage[Glenn]{fncychap}

\usepackage{pstricks}
\usepackage{pgf,tikz}
\usepackage{pstricks-add}
\usepackage{pst-plot}
\usetikzlibrary{arrows}
\usepackage{pgfplots}
\usepackage{pgfkeys}

\definecolor{qqqqqq}{rgb}{0,0,0}

\definecolor{xdxdff}{rgb}{0.49,0.49,1}
\definecolor{qqwuqq}{rgb}{0,0.39,0}
\definecolor{qqqqff}{rgb}{0,0,1}
\definecolor{ttttff}{rgb}{0.2,0.2,1}
\definecolor{uququq}{rgb}{0.25,0.25,0.25}

\renewcommand{\epsilon}{ \varepsilon}

\newcommand{\euro}{\texteuro{}}

\newcommand{\pre}{{\bf Proof.\ }}

\newcommand{\vu}{{\bf u}}

\newcommand{\vv}{{\bf v}}

\newcommand{\va}{{\bf a}}

\newcommand{\vx}{{\bf x}}

\newcommand{\cl }{\mathcal }

\newtheorem{theo}{Theorem}[section]
\newtheorem{prop}[theo]{Proposition}
\newtheorem{lem}[theo]{Lemme}
\newtheorem{cor}[theo]{Corollaire}

\newcounter{exercice}

\definecolor{fbase}{rgb}{0.8,0.8,1}
\definecolor{fgris}{gray}{0.6}
\definecolor{frouge}{HTML}{DC143C}
\definecolor{fvert}{rgb}{0.6,1,0.6}
\definecolor{fbleu}{rgb}{0.4,0.4,1}
\definecolor{fjaune}{HTML}{DCDC14}
{\endMakeFramed}
\makeatletter
\DeclareRobustCommand\sfrac[1]{\@ifnextchar/{\@sfrac{#1}}%
                                            {\@sfrac{#1}/}}
\def\@sfrac#1/#2{\leavevmode\kern.1em\raise.5ex
         \hbox{$\m@th{\fontsize\sf@size\z@\selectfont#1}$}
         \kern-.1em/\kern-.15em\lower.55ex
          \hbox{$\m@th{\fontsize\sf@size\z@\selectfont#2}$}}

\DeclareRobustCommand{\Efrac}[2]{{\displaystyle\begingroup
\raise2ex\hbox{$\m@th{#1}$}\endgroup\@@over \lower1ex
\hbox{$\m@th{#2}$}}}
\makeatother

\numberwithin{equation}{section}

\newtheorem{defi}{Definition}

\newtheorem{Exc}{Exercice}
\Newassociation{correction}{Soln}{mycor}
\Newassociation{indication}{Indi}{myind}

\def\exo#1{\futurelet\testchar\MaybeOptArgmyexoo}
\def\MaybeOptArgmyexoo{\ifx[\testchar \let\next\OptArgmyexoo
                        \else \let\next\NoOptArgmyexoo \fi \next}
\def\OptArgmyexoo[#1]{\begin{Exc}[#1]\normalfont}
\def\NoOptArgmyexoo{\begin{Exc}\normalfont}

\newcommand{\finexo}{\end{Exc}}
\newcommand{\flag}[1]{}

\tikzset{
xmin/.store in=\xmin, xmin/.default=-3, xmin=-3,
xmax/.store in=\xmax, xmax/.default=3, xmax=3,
ymin/.store in=\ymin, ymin/.default=-3, ymin=-3,
ymax/.store in=\ymax, ymax/.default=3, ymax=3,
}

\newcommand{\entete}[1]

\begin{document}
 \Opensolutionfile{mycor}[ficcorex]
 \Opensolutionfile{myind}[ficind]
 \entete{\'Enoncés}


\maketitle
\begin{abstract}
We prove that the $\mathbb{R}-$algebra $\cl S(\cl P(E,M)) $ of symbols of differential operators acting on the sections of the vector bundle $E\to M$ decompose into the sum
\[
  \cl S(\cl P(E,M))=\cl J(E)\oplus {\rm Pol}(T^*M)
\]
where $\cl J(E)$ is an ideal of $\cl S(\cl P(E,M))$ in which product of two elements is always zero.
This induces that $\cl S(\cl P(E,M))$  cannot characterize  $E \to M$ with its only structure of $\mathbb{R}-$ algebra.

We prove that with its Poisson algebra structure, $\cl S(\cl P(E,M))$ characterizes the vector bundle $E\to M$ without the requirement to be considered as a ${\rm C}^\infty(M)-$module.\\
\end{abstract}

\section{ Introduction}
An algebra $ \cl S(E)$ (Lie or associative) characterizes a vector bundle $ E $ if any isomorphism of the algebras $ \cl S(E) $ and $ \cl S(F) $ induces a diffeomorphism between the vector bundles $ E $ and $ F.$ 
When  the algebra $\cl S(E)$ is a Lie one, we speak of a Pursell-Shanks type  result; and when it is an associative algebra, one obtains a  Gel'fand-Kolmogoroff type result.\\

It is established in \cite{LecZih1} a Lie-algebraic characterization of vector bundles by the classical Poisson algebra $ \cl S(\cl P(E,M))$ of symbols of differential operators acting on the sections of a vector bundle $E\to M$. But it is required that this Lie algebra be considered as a ${\rm C}^\infty(M)-$module.\\

Classical Poisson algebras are Lie algebras but also associative algebras. The main goal in this paper is to obtain characterization of vector bundles with both Lie and $\mathbb{R}-$algebra structures of $\cl S(\cl P(E)).$ \\
The interest is that with only one of these two algebraic structures (Lie or $\mathbb{R}-$algebra) of $\cl S(\cl P(E,M))$, we don't have optimal characterization in the following sense.
\begin{enumerate}
\item[$\bullet$] The only $\mathbb{R}-$algebra structure can't characterize the vector bundle 
\item[$\bullet$] One can obtain a Lie-algebraic characterization of vector bundle with $\cl S(\cl P(E,M)),$ but this space was considered as a ${\rm C}^\infty(M)-$module. This result can be found in \cite{LecZih1}.
\end{enumerate} 
Thus, the space $\cl S(\cl P(E,M))$ offers a special opportunity  to obtain a non trivial both Pursell-Shanks and Gel'fand-Kolmogoroff type result for vector bundles.

\section{The Lie algebra of linear operators of a vector bundle}

Let $ E \to M $ be a vector bundle of rank $ n> 1$ and  denote by $ \Gamma (E) $ the space of its smooth sections.  Elements of the $\mathbb{R}-$algebra $ {\rm C}^\infty (M)$  can be seen as endomorphisms of the space $ \Gamma (E)$ by the relation  
  \[
  \gamma_u:\Gamma(E)\to\Gamma(E):s\mapsto us, \forall u\in{\rm C}^\infty(M)
  \] 
We denote the Lie algebra of differential operators by
\[
  \cl P(E,M)=\bigcup_{i\geq 0}\cl P^i(E,M),
\]
with, by definition,
\[
 \cl A(E,M)=\cl P^0(E,M)=\{\gamma_u: u\in{\rm C}^\infty(M)\}
\] 
and for $i\in\mathbb{N},$
\[
 \cl P^{i+1}(E,M)=\{T\in End(\Gamma(E))|\forall u\in{\rm C}^\infty(M): [T,\gamma_u]\in\cl P^i(E,M)\}.
\]
One can extend the filtration to $\mathbb{Z}$ by setting 
\[
\cl P^i(E,M)=\{0\},\qquad \forall i\leqslant-1.
\]
In what follows, we study the classical limit of the quantum Poisson algebra $ \cl P (E, M).$ \\
The  classical limit mentioned above is defined by
  \[
     \cl S(\cl P(E,M))=\bigoplus_{i\in\mathbb{Z}}\cl S^i(\cl P(E,M));
  \]  
with $\cl S^i(\cl P(E,M))=\cl P^i(E,M)/\cl P^{i-1}(E,M).$
We obtain a classical Poisson algebra whose operations are given in the following.
Let us begin with a definition.
\begin{defi}
Let $T\in\cl P^i(E,M).$ Then, $ord(T)=i,$ if $T\notin \cl P^{i-1}(E,M).$\\
For $i\geq ord(T),$ the $i$-degree symbol of $ T $  is defined by
\[
\sigma_i(T)=\left\{
\begin{array}{l}
 0 \mbox{ if } i>ord(T)\\
 T+\cl P^{i-1}(E,M)\mbox{ if } i=ord(T)\cdot
\end{array}
\right.
\]
\end{defi}
The structure of the quantum Poisson algebra $\cl P(E,M)$ allows to define the above usual symbol.
\begin{defi}
The symbol related to the quantum Poisson structure of $ \cl P (E, M)$ is given by
  \[
  \sigma_{pson} : \cl P(E,M)\to  \cl S(\cl P(E,M)): T\mapsto \sigma_{ord}(T).
  \]
\end{defi}
we can now precise the operations in $\cl S(\cl P(E,M)).$
\begin{defi}
Let $P\in \cl S^i(\cl P(E,M))$ and $Q\in\cl S^j(\cl P(E,M))$ such that $P=\sigma_i(T)$ and $Q=\sigma_j(D).$ We set, by definition,
 \[
  P.Q=\sigma_{i+j}(T\circ D) \mbox{ and } \{P,Q\}=\sigma_{i+j-1}([T,D])\cdot
 \]
\end{defi}

The interest of the following section is more illustrative than purely mathematical.
In fact, the mathematical considerations developed in the rest of this article have their origin in the analysis of what happens in the particular case of the algebra $gl(E).$

\section{The particular case of the $\mathbb{R}-$algebra $\sigma_{pson}(gl(E))$}
 
By virtue of the calculations made in the previous section, we have, by definition, 
 \[
 \sigma_{pson}(\gamma_u)=\gamma_u+\{0\},\quad \forall u\in{\rm C}^\infty(M)\cdot
 \]
We will then simply denote $\qquad\sigma_{pson}(\gamma_u)=\gamma_u.$\\
Likewise, for $A\in gl(E)\setminus\cl P^0(E,M),$ we have
\[
 \sigma_{pson}(A)= A'+\cl P^0(E,M),
\] 
with $A'=A-\frac{tr(A)}{n}id,$\quad $tr(A)$ being the trace of $A.$
And for the product, if $A,B\notin\cl P^0(E,M),$ we have 
  \[
  \sigma_{pson}(A)\cdot\sigma_{pson}(B)=0,
  \]
but for $\gamma_u\in\cl P^0(E,M),$ we then have   
  \[
  \sigma_{pson}(\gamma_u)\cdot\sigma_{pson}(A)=\gamma_u \circ A'+\cl P^0(E,M).
  \]
We therefore have the following identification of Lie algebras
\[
\sigma_{pson}(gl(E))\cong sl(E)\oplus {\rm C}^\infty(M)id\cdot 
\]
In the previous identification,  the bracket is given by the following relation
\begin{equation}\label{bracket sgle}
  \{A+\gamma_u , B+\gamma_v \}=[A,B].
\end{equation}
The multiplication is commutative and defined by
 \begin{equation}\label{multi sgle}
 (A+\gamma_u)\cdot(B+\gamma_v)=\gamma_v\circ A +\gamma_u\circ B +\gamma_{uv}.
 \end{equation}

\begin{prop}\label{inv sle}
When the $\mathbb{R}-$vector space $\sigma_{pson}(gl(E))=sl(E)\oplus {\rm C}^\infty(M)\,id$ is provided with the multiplication defined in (\ref{multi sgle}) above, we have
\begin{enumerate}
\item[$\bullet$]$A+\gamma_u$ is invertible if and if only $u$ is. We therefore have 
    \[
      (A+\gamma_u)^{-1}=u^{-2}(-A+\gamma_{u})
    \]
\item[$\bullet$] the subset $\cl J(E)\subset \sigma_{pson}(gl(E)) $ defined by
   \[
     \cl J(E)=\{ \va\in \sigma_{pson}(gl(E)):\va^2=0\}
   \]
 is an ideal of the $\mathbb{R}-$algebra $\sigma_{pson}(gl(E))$ and we have the decomposition 
 \[
   \sigma_{pson}(gl(E))=\cl J(E)\oplus {\rm C}^\infty(M) 
 \]  
\end{enumerate}
\end{prop}

\begin{prop}
Provided in addition with the bracket given in (\ref{bracket sgle}), the $\mathbb{R}-$algebra $\sigma_{pson}(gl(E))$ is a Poisson algebra.
\end{prop}


\section{ The $\mathbb{R}-$algebra $ \cl S(\cl P(E,M))$} 

\subsection{ From local expression}

The following is taken from \cite{LecZih1}.

\begin{prop}\label{critere pkde }
The elements of $\cl P^k(E,M),$ $k\geq1,$ are characterized by the fact that they are written locally, in a trivialization domain $ U \subset M, $ in the form
\begin{equation}\label{(*v)}
  \sum_{|\alpha|<k}T_\alpha\partial^\alpha+\sum_{|\beta|=k}u_\beta\partial^\beta
\end{equation}
where $T_\alpha\in{\rm C}^\infty (U,gl(n,\mathbb{R}))$ and $u_\beta\in{\rm C}^\infty (U).$
\end{prop}


In a trivialization domain $ U \subset M,$ the part of order strictly equal to $ k $ in the local expression of $ T\in\cl P^k(E,M)$ has the following form
\begin{equation}\label{(**)}
 \sum_{|\alpha|=k-1}A_{\alpha}\partial^\alpha+\sum_{|\beta|=k}u_{\beta}\partial^\beta, 
\end{equation}
where $A_\alpha\in {\rm C}^\infty(U,sl(n,\mathbb{R}))$  and $u_\beta\in{\rm C}^\infty(U).$\\
Let now consider $D\in\cl P^l(E,M), $ whose $l-$order terms are given by
\[
 \sum_{|\lambda|=l-1}B_{\lambda}\partial^\lambda+\sum_{|\mu|=l}v_\mu\partial^\mu  \]
Therefore, the terms of exactly $(k+l)-$order of $T\circ D$  are given by 
\begin{equation}\label{produit loc}
   \sum_{|\alpha|=k-1}\sum_{|\mu|=l}A_{\alpha}v_\mu\partial^\alpha\partial^\mu
   +
  \sum_{|\lambda|=l-1}\sum_{|\beta|=k} B_\lambda u_\beta\partial^\beta\partial^\lambda
   +
    \sum_{|\beta|=k}\sum_{|\mu|=l}u_\beta v_\mu\partial^\beta\partial^\mu
\end{equation}

We notice that the local decomposition (\ref{(**)})  given above is not intrinsic.
In fact, if in the sum on the right we recognize the principal symbol of the differential operator in the usual sense, the sum on the left, for its part, does not resist a change of coordinates and is therefore not globally defined.

In the following lines we build a global decomposition allowing to find a global meaning to the expression given in (\ref{(**)})  previously.\\

Let us first recall that the Lie algebra of linear operators of the vector bundle $E\to M$ can be defined by its usual filtration as follows:
  \[
    \cl D^0(E,M)=\{T\in End(\Gamma(E)):[T,\gamma_u]=0,\,\forall u\in{\rm C}^\infty(M)\}
  \]
For any integer $k\geq1,$
  \[  
    \cl D^k(E,M)= \{T\in End(\Gamma(E))| \forall u\in{\rm C}^\infty(M): [T,\gamma_u]\in    
    \cl D^{k-1}(E,M)\}.
  \]
We then set
 \[
   \cl D(E,M) =\bigcup_{k\geq 0} \cl D^k(E,M).
 \]
Let now $T\in\cl S^{k-1}(M)\otimes sl(E)$ and assume given a unit partition $(U_i,\rho_i)$ of $M$ whose domains $U_i$ are the trivialization one of $E.$ In any $U_i,$ assume $T$ is expressed in the form
 \[
   T=\sum_{|\alpha|=k-1}A_{\alpha,i}\xi^\alpha\cdot
 \]
We therefore set
 \[
  \overline{T}_i=\sum_{|\alpha|=k-1}A_{\alpha,i}\partial^\alpha\in\cl D^{k-1}(E|_{U_i},U_i)
 \]
with $A_{\alpha,i}\in  {\rm  C}^\infty(U_i,sl(n,\mathbb{R})).$ 
We obtain the following differential operator
  \[
    \overline{T}=\sum_i \rho_i\overline{T}_i\in \cl D^{k-1}(E,M)\subset\cl P^k(E,M).
  \]
It is associated with the partition of the unit chosen at the start and it is then such that
   \[
    \sigma_{pson}(\overline{T})=\sigma_{ppal}(\overline{T})=T, \footnote{Note that $\sigma_{ppal}$ is the usual principal symbol of differential operators.}
   \]
but it is obviously not the only one to verify this relation. \\
Nevertheless, we have the following statement.
\begin{prop}
The space $\cl S^k(\cl P(E,M))=\cl P^k(E,M)/\cl P^{k-1}(E,M)$ of symbols in the sense "quantum Poisson algebra" of differential operators in $\cl P^k(E,M)$ is determined by the following  exact short sequence of $\mathbb{R}-$vector spaces 
   \[
    0\longrightarrow\cl S^{k-1}(M)\otimes sl(E)\stackrel{\theta}{\longrightarrow}\cl P^k(E,M)/\cl P^{k-1}(E,M)\stackrel{\delta}{\longrightarrow}\cl S^k(M)\longrightarrow 0,
   \] 
with $\theta: T\mapsto \overline{T}+\cl P^{k-1}(E,M)$ and 
  \[
  \delta: D+\cl P^{k-1}(E,M)\mapsto \left\{
                                 \begin{array}{l}
                                 0 \mbox{ if } D\in\cl D^{k-1}(E,M)\\
                                 \sigma_{ppal}(D) \mbox{ if not. }
                                 \end{array}
                                \right. 
  \] 
\end{prop}
Therefore, seen as $\mathbb{R}-$vector spaces, we have the following decomposition 
  \[
    \cl S^k(\cl P(E,M))=Pol^{k-1}(T^*M,sl(E))\oplus Pol^k(T^*M,\mathbb{R})
  \]
for any integer   $k\in\mathbb{N}.$

\subsection{ Global point of view}

For a vector bundle $ E \to M,$ of rank $n\geqslant 2$, we have established in \cite{LecZih1} that under certain assumptions, the Lie structure of $\cl S(\cl P(E,M))$ characterizes the vector bundle $ E.$ We will see what it is about its structure of $\mathbb{R}-$ algebra.\\

Let us consider the following exact sequence, encountered in the previous section, and let's focus to the only associative algebra structure of the spaces in question:
\[
0\to \cl S(M)\otimes sl(E)\stackrel{\theta}{\longrightarrow}\cl S(\cl P(E,M))\stackrel{\delta}{\longrightarrow} \cl S(M)\to 0.
\]
Note that by virtue of the relation (\ref{produit loc}) in the previous section, the structure of $\mathbb{R}-$ algebra of $\cl S(M) \otimes sl(E)$ in question here is the one for which the product of any two elements is zero.
 
Let us consider the subset of $\cl S(\cl P(E,M))$ defined and denoted by
  \[
     \cl J(E):=\{P\in\cl S(\cl P(E,M)):P^2=0\}\cdot
  \] 
We then have that the space $\cl J(E)$ is a ideal of the $\mathbb{R}-$algebra $\cl S(\cl P(E,M)).$
Indeed, this follows from the fact that $ \cl J(E) $ is nothing other than the kernel of the homomorphism of $ \mathbb{R}-$ algebras $ \delta$ given in the preceding exact sequence.\\
As in the particular case of the algebra $ \sigma_{pson}(gl(E)),$ let us determine the invertible elements of the $\mathbb{R}-$algebra $\cl S(\cl P(E,M)).$
\begin{prop} 
The invertible elements of the $ \mathbb{R}-$ associative algebra $\cl S(\cl P(E,M))$ decompose in the form $ \vu + f, $ with $ \vu \in \cl J(E) $ and $ f \in {\rm Pol}^0 (T^* M) \cong {\rm C}^\infty(M)$ a nonvanishing function. And the inverse of such an element is given by
   \[ 
     (\vu+f)^{-1}=-f^{-2}\cdot\vu+f^{-1}.
   \]  
with $f^{-1}:T^*M\to \mathbb{R}: \vx\mapsto \frac{1}{f(\vx)}.$
\end{prop}
\pre
  Let  $P,Q\in \cl S(\cl P(E,M))$ such that $P\cdot Q=1.$ 
  We then have 
  \[
  \delta(P)\cdot\delta(Q)=1,
  \]
and this implies that $ \delta (P) $ and $ \delta (Q) $ are constant polynomials along the fibers of $ T^*M, $ taking into account the identification $ \cl S (M) \cong {\rm Pol} (T^* M). $\\ We deduce that there are non vanishing functions $ f, g \in {\rm C}^\infty(M),$ one being the inverse of the other , such as $\delta(P)=f$ and $\delta(Q)=g.$ \\
Therefore, by definition of the homomorphism $\delta,$ we have the decomposition
\[
  P=f+\vu \mbox{ and } Q=g+\vv
\]
with $\vu,\vv\in\cl J(E)=ker\,\delta.$
We then obtain
  \begin{eqnarray*}
    (f+\vu)^{-1} & = & \frac{1}{f}\frac{1}{1+\frac{1}{f}\vu}\\
               & = & \frac{1}{f}\sum_k (-1)^k\frac{1}{f^k} \vu^k= f^{-1}-f^{-2}\vu
  \end{eqnarray*}  
since $\vu^2=0,$ 
and the result is established. \hfill $\blacksquare$ \\
This result generalizes what we have obtained in the Proposition \ref{inv sle}.\\

Observe that the gradation of the $ \mathbb{R}-$ algebra $ \cl S(\cl P (E, M)) $ induces one on $ \cl J(E) $ and we can write
\[
\cl J(E)=\bigoplus_{k\geq 0}\cl J^k(E),
\]
with in particular $\cl J^0(E)=\{0\}.$\\

We can now state the following result.
\begin{prop}\label{ ideal de spe}
Let $ E \to M $ and $ F \to N $ be two vector bundles. If $ \Psi: \cl S (\cl P (E, M)) \to \cl S (\cl P (F, N)) $ is an isomorphism of $ \mathbb{R}-$ algebras, then $ \Psi $ respects their ideals $\cl J(E)$ and $\cl J(F)$. 
\end{prop} 
\pre
Let $\Psi:\cl  S(\cl P(E,M))\to\cl S(\cl P(F,N))$ be an isomorphism of $\mathbb{R}-$algebras. 
The proposition comes directly from the fact that for any $ P \in \cl S (\cl P (E, M)) $ such that there exists $ r \in \mathbb{N} $ with $ P^r = 0,$ we also have 
$(\Psi(P))^r=0.$\hfill $\blacksquare$ \\

The result below is taken from \cite{Zih}.
\begin{lem}\label{lem:degree 0}
Let $E \to M$ and $F \to N$ be two vector bundles.\\
If $\Psi: {\rm Pol}(E) \to {\rm Pol}(F)$ is an isomorphism of associative algebras, we then have
 \[
  \Psi({\rm Pol}^0(E)) = {\rm Pol}^0(F)
 \]
\end{lem}
 
We can now state the following result which induces a sort of algebraic-characterization of manifolds.
 
\begin{prop}\label{identification polm et poln}
 
Let $ E \to M $ and $ F \to N $ be two vector bundles. Any isomorphism of $ \mathbb{R}-$algebras $\Psi:\cl S(\cl P(E,M))\to \cl S(\cl P(F,N))$  satisfies the following relation
   \[
  \Psi(\cl S^0(\cl P(E,M)))=\cl S^0(\cl P(F,N)).
  \]  
\end{prop} 
\pre
Since $\cl J(E)$ (resp. $\cl J(F)$) is an ideal of the $\mathbb{R}-$algebra  of $\cl S(\cl P(E,M))$ (resp. $\cl S(\cl P(F,N))$) and $\Psi(\cl J(E))=\cl J(F),$ we then deduce that
\[
  \overline{\Psi}: [Q]\mapsto [\Psi(Q)]
\]
is indeed an isomorphism of $ \mathbb{R}-$ algebras between the quotient spaces\\ 
$ \cl S(\cl P(E,M))/\cl J(E)$ and $\cl S(\cl P(F,N))/\cl J(F).$ \\
Observe that these quotient $\mathbb{R}-$algebras  are graded, the subspace of $k-$weight  being identified with
  \[
    \cl S^k(\cl P(E,M))/\cl J^k(E)\cong {\rm Pol}^k(T^*M)\cdot
  \]
For  $\cl S(\cl P(F,N))/\cl J(F),$ we have a similar identification.
 We deduce the existence of graded isomorphisms of $\mathbb{R}-$algebras
  \[
    \overline{\Psi}_M: \cl S(\cl P(E,M))/\cl J(E)\to {\rm Pol}(T^*M)
  \]
   and  
  \[
  \overline{\Psi}_N: \cl S(\cl P(F,N))/\cl J(F)\to {\rm Pol}(T^*N).  
  \] 
We then get an isomorphism of $\mathbb{R}-$algebras
  \[
    \overline{\Psi}_N\circ\overline{\Psi}\circ\overline{\Psi}_M^{-1}: {\rm Pol}(T^*M)\to {\rm Pol}(T^*N).
  \]  
By virtue of the previous Lemma \ref{lem:degree 0}, we have
\[
\overline{\Psi}_N\circ\overline{\Psi}\circ\overline{\Psi}_M^{-1}({\rm Pol}^0(T^*M))={\rm Pol}^0(T^*N).
\]
This allows to write
   \[
     \overline{\Psi}( \cl S^0(\cl P(E,M))/\cl J^0(E))=\cl S^0(\cl P(F,N))/\cl J^0(F),
   \]
since $\overline{\Psi}_M$ and $\overline{\Psi}_N$ are graded. 
Therefore, for any $f\in\cl S^0(\cl P(E,M)),$ we have the equality
  \[
    \overline{\Psi} ([f])=[\Psi(f)]\in\cl S^0(\cl P(F,N))/\{0\},
  \]  
and the desired result follows. \hfill $\blacksquare$ \\

We deduce the following result.
\begin{cor}
  Let $E\to M$ and $F\to N$ be two vector bundles. 
If the $ \mathbb{R}-$ algebras $ \cl S(\cl P(E,M)) $ and $ \cl S(\cl P(F,N)) $ are isomorphic, then the manifolds $ M $ and $ N $ are diffeomorphic. \\
\end{cor}

Observe that we can't do more. Indeed, a vector bundle $E \to M$ cannot be characterized by the only associative algebra structure of $\cl S(\cl P(E,M))$ since, 
"the information" on the nature of the bundle, except its base, is housed in the ideal $ \cl J(E) $ on which the restriction of the multiplication of $\cl S(\cl P(E,M)) $ reduces to the trivial structure.


\section{ Poisson-algebraic characterization}

For quantum Poisson algebras, a homomorphism of $ \mathbb{R}-$ algebras is necessarily a homomorphism of Lie algebras. \\
The results obtained in the previous section and in \cite{GraPon4,LecZih,LecZih1,Zih} suggest to consider, for classical Poisson algebras, homomorphisms respecting both associative and Lie algebra structures, in order to obtain an algebraic characterization of vector bundles.\\
As we show it in the following lines, we then don't need to impose that the $ \mathbb {R}-$ algebras $ \cl S (\cl P(E, M)) $ and $ \cl S (\cl P(F, M)) $ are seen as ${\rm C}^\infty(M)-$modules.\\

But first let's give a definition.

\begin{defi}
A $ \mathbb{R}-$ linear map $ \Psi $  between two classical Poisson algebras is a \textsl{homomorphism of Poisson algebras} if $ \Psi $ is both a homomorphism of associative algebras and of Lie algebras. 
\end{defi}

Here, $\Psi$ is not required to be graded.\\

\begin{lem}\label{filtr}
Denoting $\cl S(\cl P(F,N))$ by $\cl S,$ we have
\begin{equation*}
   \{T\in\cl S:\{T,f\}\in \cl S^{i-1},\forall f\in\cl A\}=\cl S^i\oplus gl(F), \forall i\in\mathbb{N}.
\end{equation*}
where $\cl A$ is $\cl A(F,N).$
\end{lem}
\pre
The inclusion 
        \[\{P\in\cl S:\{P,\cl A\}\subset\cl S^{i-1}\}\supset\cl S^{i}\oplus gl(F) \]
being immediate, it remains to establish reciprocal inclusion.
Let $ P\in\cl S$ such that $\{P,\cl A\}\subset \cl S^{i-1}$ and assume absurdly that 
     \[P\notin gl(F)\oplus\cl S^{i}.\]
We can then write  
     \[P= P_0+\Sigma_k H_k\] 
where $P_0\in gl(F)\oplus\cl S^{i}$, $H_k\in\cl S^k-\{0\},k\neq i,$ and $H_k\notin gl(F).$
We deduce, on the one hand, that 
  \[\{H_k,f\}\in\cl S^{k-1},\forall f\in\cl A \cdot\]
On the other hand, for all $f\in\cl A$, we obtain
   \[\Sigma_k\{H_k,f\}\in\cl S^{i-1},\] 
since $\{P,f\}-\{P_0,f\}=\Sigma_k\{H_k,f\}.$ 
We therefore have 
     \[\{H_k,f\}=0, \forall f\in\cl A.\]
This implies that $H_k\in gl(F).$ It is absurd. \hfill $\blacksquare$\\

We can now state the following proposition.
\begin{prop}\label{iso filtr}
Let $E\to M$ and $F\to N$ be two vector bundles.\\
Any isomorphism of Poisson algebras between $\cl S (\cl P(E, M))$ and $\cl S(\cl P(F,N))$ respects the gradation.
\end{prop}   
\pre
Let $\Psi:\cl S(\cl P(E,M))\to \cl S(\cl P(F,N))$ be an isomorphism of Poisson algebras.
Observe that by virtue of the previous Proposition \ref{identification polm et poln}, we have
  \begin{equation}\label{resp filtr 0}
  \Psi (\cl S^0(\cl P(E,M)))=\cl S^0(\cl P(F,N)).
  \end{equation}
Let us now assume, by the induction hypothesis, that
 \[\Psi(\cl S^i(\cl P(E,M)))\subset\cl S^i(\cl P(E,N)).\] For any $T\in\cl S^{i+1}(\cl P(E,M))$ we have $\{T,\gamma_u\}\in\cl S^i(\cl P(E,M)),$ for all $u\in{\rm C}^\infty(M).$
This implies that \[\{\Psi(T),\gamma_v\}\in\cl S^i(\cl P(E,N)),\] for any $v\in{\rm C}^\infty(N),$ according to the induction hypothesis and to the equality (\ref{resp filtr 0}).
We therefore obtain, by using the Lemma \ref{filtr}, 
\[\Psi(T)\in\cl S^{i+1}(\cl P(E,N)).\]
This conclude the proof of the proposition. \hfill $\blacksquare$ \\

The following lemma is taken from \cite{Lec3}. 

\begin{lem}\label{ende carct e}
Let $ E \mapsto M $ and $ F \mapsto M $ be two vector bundles of respective ranks $n,n'>1$ with $H^1(M,\mathbb{Z}/2)=0$. The Lie algebras $ gl (E) $ and $ gl (F) $ (resp. $ sl (E) $ and $ sl (F) $) are isomorphic if and only if the vector bundles $ E $ and $ F $ are isomorphic.
\end{lem}

The following proposition gives a Poisson-algebraic characterization of vector bundles.

\begin{theo}
Let $E\to M$ and $F\to N$ be two vector bundles of respective ranks $n,n'>1$.\\
Under the hypothesis $ H^1(M,\mathbb{Z}/2) = 0, $ the Poisson algebras $\cl S (\cl P(E,M))$ and $\cl S(\cl P(F,N))$ are isomorphic if, and only if, the vector bundles $ E $ and $ F $ are.
\end{theo}   
\pre
Let $\Psi:\cl S(\cl P(E,M))\to \cl S(\cl P(F,N))$ be an isomorphism of Poisson algebras. 
Let us show that
\[
  \Psi(sl(E))=sl(F).
\] 
This is sufficient to conclude according to the Lemma \ref{ende carct e}.\\ 
For any $A\in sl(E)\subset \cl S^1(\cl P(E,M)),$ the Proposition \ref{iso filtr} allows us to write $\Psi(A)\in \cl S^1(\cl P(F,N))$. We also have, 
\begin{equation*}\label{resp filtr 1}
 \{\Psi(A),\gamma_v\}=0, \forall v\in{\rm C}^\infty(N),
\end{equation*} 
by virtue of the previous relation (\ref{resp filtr 0}). This implies $\Psi(A)\in sl(F).$ \hfill $\blacksquare$ \\



\end{document}